\documentclass[11pt]{amsart}  
\usepackage{amsmath,amssymb,amsfonts,times,hyperref}

\title[Codes in spherical caps]{Semidefinite programming, multivariate
orthogonal polynomials, and codes in spherical caps}

\author{Christine Bachoc} 
\address{C. Bachoc, Laboratoire A2X, Universit\'e Bordeaux I, 351,
cours de la Li\-b\'e\-ration, 33405 Talence, France}
\email{bachoc@math.u-bordeaux1.fr}

\author{Frank Vallentin} 
\address{F. Vallentin, Centrum voor Wiskunde en Informatica (CWI),
Kruislaan 413, 1098 SJ Amsterdam, The Netherlands}
\email{f.vallentin@cwi.nl}

\thanks{The second author was supported by the Netherlands
Organization for Scientific Research under grant NWO 639.032.203 and
by the Deutsche Forschungsgemeinschaft (DFG) under grant SCHU
1503/4.}


\date{October 16, 2007}
\newtheorem{defi}{Definition}[section]
\newtheorem{definition}[defi]{Definition}
\newtheorem{proposition}[defi]{Proposition}
\newtheorem{theorem}[defi]{Theorem}
\newtheorem{remark}[defi]{Remark}
\newtheorem{corollary}[defi]{Corollary}
\newtheorem{lemma}[defi]{Lemma}

\newcommand{\R}{{\mathbb{R}}}

\newcommand{\I}{{\mathcal{I}}} 
\newcommand{\mymid}{:} 

\newcommand{\Sn}{S^{n-1}} 
 
\newcommand{\Pd}{\operatorname{Pol}_{\leq d}}

\newcommand{\Pl}{\operatorname{Pol}}
\newcommand{\On}{{\operatorname{O}(\R^n)}}

\newcommand{\Trace}{\operatorname{Trace}}
\newcommand{\Harm}{\operatorname{Harm}}

\newcommand{\card}{\operatorname{card}}
\newcommand{\Stab}{\operatorname{Stab}}

\newcommand{\vol}{\operatorname{vol}} 

\newcommand{\prodeucl}[2]{#1 \cdot #2}
\newcommand{\prodhaar}[2]{(#1,#2)}

\newcommand{\CP}{\operatorname{Cap}(e,\phi)}
\newcommand{\Y}{\overline{Y}_k^n}

\begin{document}

\begin{abstract} 
In this paper we apply the semidefinite programming approach developed
in \cite{BV} to obtain new upper bounds for codes in spherical
caps. We compute new upper bounds for the one-sided kissing number in
several dimensions where we in particular get a new tight bound in
dimension $8$. Furthermore we show how to use the SDP framework to get
analytic bounds.
\end{abstract}

\maketitle

\centerline{\large \em Dedicated to Eiichi Bannai in occasion of his
  60th birthday}
\section{Introduction}

Let $S^{n-1}$ denote the unit sphere of the Euclidean space
$\R^n$. The spherical cap with center $e\in S^{n-1}$ and angular
radius $\phi$ is the set 
\begin{equation*}
\CP=\{x\in S^{n-1} : e\cdot x \geq \cos\phi\}.
\end{equation*}
We consider the problem of finding upper bounds of the size of a code
$C$ contained in $\CP$ with minimal angular distance $\theta$.
Following notations of \cite{BM}, the maximal size of such a code is
denoted by $A(n,\theta,\phi)$.  Many reasons to consider this problem
are exposed in \cite{BM}, e.g.\ upper bounds for spherical codes can
be derived from upper bounds for spherical cap codes through the
following inequality:
\begin{equation*}
\frac{A(n,\theta)}{\vol(S^{n-1})}\leq \frac{A(n,\theta,\phi)}{\vol(\CP)}
\end{equation*}
where $A(n,\theta)$ stands as usual for the maximal size of a
spherical code with minimal angular distance $\theta$.

Moreover, it is a challenging problem, because the so-called linear
programming method does not apply to this situation.  In coding theory
many of the best upper bounds are consequences of the so-called linear
programming method due to P.~Delsarte. This method gives upper bounds
for codes from the solution of a certain linear program. It can be
applied to symmetric spaces and has been successfully used to deal
with two-point homogeneous spaces like the unit sphere $S^{n-1}$
(\cite{D}, \cite{DGS}, \cite{KL} and the survey \cite[Chapter 9]{CS}), or with symmetric spaces
which are not two-point homogeneous like Grassmannian spaces
(\cite{B}). However the method is not applicable to spaces which are
not symmetric spaces like spherical caps.

In this paper, we show that the approach developed in
\cite{BV} based on semidefinite programming can be applied to the
above problem. It turns out that it gives good numerical results. In
particular we obtain improvements in the determination of the
so-called one-sided kissing number, corresponding to $\phi=\pi/2$ and
$\theta=\pi/3$, and denoted by $B(n)$ after \cite{M2}.

Let us describe briefly the idea underlying our approach.  The
isometry group of $\CP$ is the group $H:=\mathrm{Stab}(\On,e)$
stabilizing the point $e$ in $\On$.  This group acts on the space
$\Pd(S^{n-1})$ of polynomial functions on the unit sphere of degree at
most $d$. In the decomposition of this space into irreducible
subspaces some irreducible subspaces occur with multiplicities.  To
each irreducible subspace with multiplicity $m$ we can associate an
$m\times m$ matrix $Y$ whose coefficients are real polynomials in
three variables $(u,v,t)$ and have an explicit expression in terms
of Gegenbauer polynomials. Each matrix $Y$ satisfies the positivity
property:

\begin{equation*}
\text{For all finite $C\subset \Sn$,} \sum_{(c,c')\in
  C^2}Y(\prodeucl{e}{c},\prodeucl{e}{c'}, \prodeucl{c}{c'}) \succeq 0,
\end{equation*}

\noindent
where ``$\succeq 0$'' stands for ``is positive semidefinite''.

We want to point out that one can consider other metric spaces $X$
with isometry group in this framework. Only the expression of
the matrices $Y$ will depend on the specific situation. For a
symmetric space $X$ the multiplicities in the irreducible
decomposition are equal to $1$. Hence the matrices $Y$ have size $1
\times 1$. So we recover the classical positivity property of zonal
polynomials which underlies the linear programming method.

The paper is organized as follows: Section~\ref{review} recalls the
needed notations and results of \cite{BV}. Section~\ref{sdp} states
the semidefinite program (SDP for short) which obtains an upper bound
for $A(n, \theta, \phi)$ and presents numerical results.
Section~\ref{polynomials} translates the dual SDP into a statement on
three variable polynomials, and contains more material on
orthogonality relations, positivity property and other classical
material which might be of independent interest.

\section{Review on the semidefinite zonal matrices}
\label{review}

We start with some notations. The standard inner product of
the Euclidean space $\R^n$ is denoted by $\prodeucl{x}{y}$. The
orthogonal group $\On$ acts homogeneously on the  unit
sphere
\[
\Sn := \{x\in \R^n \mymid \prodeucl{x}{x}=1\}.
\]
The space of real polynomial functions on $\Sn$ of degree at
most $d$ is denoted by $\Pd(\Sn)$. It is endowed with the
induced action of $\On$, and equipped with the standard
$\On$-invariant inner product
\begin{equation}\label{ip}
\prodhaar{f}{g}=\frac{1}{\omega_n}\int_{\Sn} f(x)g(x)d\omega_n(x),
\end{equation}
where $\omega_n$ is the surface area of $\Sn$ for the standard measure
$d\omega_n$.

It is a classical result that under the action of $\On$
\begin{equation}\label{dec 1}
\Pd(\Sn)=H_0^n\perp H_1^n\perp\ldots\perp H_d^n,
\end{equation}
where $H_k^n$ is isomorphic to the $\On$-irreducible space of
homogeneous, harmonic polynomials of degree $k$ in $n$ variables,
denoted by $\Harm_k^n$.  For the dimension of these spaces we write
$h_k^n:=\dim(\Harm_k^n)$.

For the restricted action of the subgroup $H:=\Stab(e,\On)$,
introduced above, we have the following decomposition into isotypic
components:
\begin{equation}\label{dec iso}
\Pd(\Sn)=\I_0\perp \I_1 \perp \ldots \perp \I_d,
\end{equation}
where
\begin{equation*}
\I_k\simeq
(d-k+1)\Harm_k^{n-1}, \quad \mbox{$k = 0, \ldots, d$.}
\end{equation*}
More precisely, $\I_k$ decomposes as 
\begin{equation}\label{dec Ik}
\I_k=H_{k,k}^{n-1} \perp \ldots \perp H_{k,d}^{n-1},
\end{equation}
where, for $i\geq k$, $H_{k,i}^{n-1}$ is the unique subspace of
$H_i^n$ isomorphic to $\Harm_k^{n-1}$.

The following construction associates to each $\I_k$ a matrix-valued
function
\begin{equation}
Z^n_k : \Sn \times \Sn \to \R^{(d-k+1)\times (d-k+1)}
\end{equation}
which is uniquely defined up to congruence. Let $(e_{s,1}^k,
e_{s,2}^k,\dots, e_{s,h_{k}^{n-1}}^k)$ be an orthonormal basis of
$H_{k,k+s}^{n-1}$. We assume that the basis $(e_{s,i}^k)_{1\leq i\leq
h_{k}^{n-1}}$ is the image of $(e_{0,i}^k)_{1\leq i\leq h_{k}^{n-1}}$
by some $H$-isomorphism $\phi_s:H_{k,k}^{n-1} \to
H_{k,k+s}^{n-1}$. Then, define
\[E_k^n(x):=
\frac{1}{\sqrt{h_k^{n-1}}}\begin{pmatrix}
e_{0,1}^k(x) & \ldots & e_{0,h_k^{n-1}}^k(x)\\
\vdots&&\vdots\\
e_{d-k,1}^k(x) & \ldots & e_{d-k,h_k^{n-1}}^k(x)
\end{pmatrix},
\]
and
\begin{equation}\label{Z}
Z_k^n(x,y):=E_k^n(x)E_k^n(y)^t \in \R^{(d-k+1) \times (d-k+1)}.
\end{equation}

One can prove that, for all $g\in H$, $Z_k^n(g(x),g(y))=Z_k^n(x,y)$.
As a consequence, the coefficients of $Z_k^n$ can be expressed as
polynomials in the three variables $u=e\cdot x$, $v=e\cdot y$,
$t=x\cdot y$.  More precisely, let $Y_k^n(u,v,t)$ be the
$(d-k+1)\times (d-k+1)$ matrix such that
\begin{equation}\label{add form 2}
Z_k^n(x,y)=Y_k^n(\prodeucl{e}{x},\prodeucl{e}{y}, \prodeucl{x}{y}).
\end{equation}

We denote the zonal polynomials of the unit sphere $S^{n-1}$ by
$P_k^n$. In other words, if $n \geq 3$, then $P_k^n(t)$ is the
Gegenbauer polynomial of degree $k$ with parameter $n/2-1$, normalized
by the condition $P_k^n(1)=1$. If $n = 2$, then $P_k^n(t)$ is the
Chebyshev polynomial of the first kind with degree $k$. We give in
\cite[Theorem 3.2]{BV} the following explicit expressions for the
coefficients of the matrices $Y_k^n$:

\begin{theorem}\label{Th Y}  We have, for all $0\leq i,j\leq d-k$,
\begin{equation}\label{Y}
\big(Y_k^n\big)_{i,j}(u,v,t)= \lambda_{i,j} P_{i}^{n+2k}(u)P_{j}^{n+2k}(v)Q_k^{n-1}(u,v,t),
\end{equation}
where 
\[
\displaystyle
Q_k^{n-1}(u,v,t):=\big((1-u^2)(1-v^2)\big)^{k/2}P_k^{n-1}\Big(\frac{t-uv}{\sqrt{(1-u^2)(1-v^2)}}\Big),
\]
and
\[
\lambda_{i,j}=\frac{\omega_n}{\omega_{n-1}}\frac{\omega_{n+2k-1}}{\omega_{n+2k}}(h_i^{n+2k}h_j^{n+2k})^{1/2}.
\]
\end{theorem}

We recall the matrix-type positivity property of the matrices $Y_k^n$
which underlies the semidefinite programming method:

\begin{theorem}\label{Th pos Y}
For any finite code $C\subset S^{n-1}$,
\begin{equation}\label{pos Y}
\sum_{(c,c')\in
  C^2}Y_k^n(\prodeucl{e}{c},\prodeucl{e}{c'}, \prodeucl{c}{c'})\succeq 0.
\end{equation}
\end{theorem}

\begin{proof} We recall the straightforward argument:
\[\displaystyle \sum_{(c,c')\in C^2}Z_k^n(c,c')=\Big(\sum_{c\in
C}E_k^n(c)\Big)\Big(\sum_{c\in C}E_k^n(c)\Big)^t\succeq 0.\]
\end{proof}

\section{Semidefinite programming bound for codes in spherical caps}
\label{sdp}

Let $C\subset \CP$ be a code of minimal angular distance
$\theta$. Define the domains $\Delta$ and $\Delta_0$ by
\[
\begin{array} {ll}
\Delta := \{(u,v,t) : &\cos\phi \leq u \leq v\leq 1, \\[1ex]
                   &-1\leq t\leq \cos\theta, \\[1ex]
                   &1+2uvt-u^2-v^2-t^2 \geq 0\},
\end{array}
\]
and 
\[
\Delta_0 := \{(u,u,1) : \cos\phi \leq u\leq 1\}.
\]
The two-point distance distribution of $C$ is the map $y : \Delta \cup
\Delta_0 \to \R$ given by
\begin{equation*}
y(u,v,t)=\frac{m(u,v)}{\card(C)}\card\{(c,c')\in C^2: \prodeucl{e}{c}=u,
\prodeucl{e}{c'}=v, \prodeucl{c}{c'}=t\},
\end{equation*}
where
\[m(u,v)=\begin{cases}
2 \quad \text{ if } u \neq v,\\
1 \quad \text{ if } u = v.
\end{cases}
\]
We introduce the symmetric matrices $\Y(u,v,t)$ defined by
\begin{equation*}
\Y(u,v,t):=\frac{1}{2}\Big(Y_k^n(u,v,t)+Y_k^n(v,u,t)\Big).
\end{equation*}
Then, \eqref{pos Y} is equivalent to the semidefinite condition
\[
\sum_{(u,v,t) \in \Delta \cup \Delta_0} y(u,v,t) \Y(u,v,t) \succeq 0.
\]
For any $d \geq 0$, the $y(u,v,t)$'s satisfy the following obvious
properties:
\[
\begin{array}{ll}
& \displaystyle y(u,v,t)\geq 0 \text{ for all } (u,v,t) \in \Delta \cup \Delta_0,\\[1ex]
& \displaystyle y(u,v,t) = 0 \text{ for all but finitely many } (u,v,t) \in \Delta \cup \Delta_0,\\[1ex]
& \displaystyle\sum_{(u,u,1)\in \Delta_0} y(u,u,1)=1,\\
& \displaystyle\sum_{(u,v,t)\in \Delta\cup \Delta_0} y(u,v,t)=\card(C),\\
& \displaystyle \sum_{(u,v,t)\in \Delta\cup \Delta_0}  y(u,v,t) \Y(u,v,t)
  \succeq 0
\text{ for } k = 0, \ldots d.
\end{array}
\]
Hence a solution to the following semidefinite program is an upper
bound for $A(n,\theta, \phi)$.
\[
\begin{array}{ll}
\max\big\{
& \displaystyle 1 + \sum_{(u,v,t)\in \Delta} y(u,v,t) \mymid\\[2ex]
& \mbox{$\displaystyle y(u,v,t)\geq 0\;$ for all $(u,v,t)\in \Delta\cup \Delta_0$,}\\[1ex]
& \mbox{$\displaystyle y(u,v,t)=0\;$  for all but finitely many $(u,v,t)\in \Delta\cup \Delta_0$,}\\[1ex]
& \mbox{$\displaystyle \sum_{(u,u,1)\in \Delta_0} y(u,u,1)=1$,} \\
& \mbox{$\displaystyle \sum_{(u,v,t)\in \Delta\cup \Delta_0} y(u,v,t)\Y(u,v,t)\succeq 0\;$ for all $k = 0, \ldots, d$}
\big\}.
\end{array}
\]

\medskip

As usual, the dual problem is easier to handle. The duality theorem
says that any feasible solution of the dual problem provides an upper
bound for $A(n,\theta,\phi)$. For expressing the dual problem we use
the standard notation $\langle A, B\rangle = \Trace(AB^t)$.

\begin{theorem}\label{Th SDP} Any feasible solution to the following
  semidefinite problem provides an upper bound on $A(n,\theta,\phi)$.
\begin{equation}
\label{dual sdp}
\begin{array}{ll}
\min\big\{ & 1 +M \mymid\\[1ex]
& \mbox{$F_k \succeq 0\;$ for  all $k = 0, \ldots, d$,}\\[1ex]
& \displaystyle \sum_{k=0}^d \langle F_k,
\Y(u,u,1) \rangle\leq M\;\; \text{ for all } (u,u,1)\in \Delta_0,\\[1ex]
& \displaystyle \sum_{k=0}^d \langle F_k, \Y(u, v,t)\rangle \leq -1
\;\; \text{ for all } (u,v,t)\in \Delta\big\}
\end{array} 
\end{equation}
\end{theorem}

In order to make use of this theorem in computations we follow the
same line as in \cite[Section 5]{BV}. A theorem of M.~Putinar (\cite{Pu})
shows that the two last conditions can be replaced by:
\[
\begin{array}{ll}
& \displaystyle \sum_{k=0}^d \langle F_k,
\Y(u,u,1)= M - q_0(u)-p(u)q_1(u)\\[1ex]
& \displaystyle \sum_{k=0}^d \langle F_k, \Y(u, v,t)\rangle = -1
-r_0(u,v,t)-\sum_{i=1}^4 p_i(u,v,t)r_i(u,v,t)
\end{array} 
\]
where $p(u)= -(u-\cos\phi)(u-1)$, $p_1=p(u)$, $p_2=p(v)$, $p_3=
-(t+1)(t-\cos\theta)$, $p_4 = -(u^2+v^2+t^2) + 2uvt + 1$, and the
polynomials $q_i(u)$, $0\leq i\leq 1$ and $r_i(u,v,t)$, $0\leq i\leq
4$ are sums of squares of polynomials. If we set the degree of those
polynomials to be less than a given value $N$, and fix the parameter
$d$, we relax \eqref{dual sdp} to a finite semidefinite program.

In the most interesting case $\cos\phi=0$ and $\cos\theta=1/2$,
corresponding to the so-called one-sided kissing number $B(n)$, we
obtain the computational results given in Table 1. For our
computations we chose the parameter $d = N = 10$.

\begin{table}
\label{Table kissing}
\begin{tabular}{c|c|l|c}
    & best lower  & best upper bound & SDP\\
$n$ & bound known & previously known &  method \\
\hline
3  & 9 & 9 \cite{F} & 9 \\
4  & 18 & 18 \cite{M2} & 18 \\
5  & 32 & 35 \cite{M3} & 33\\
6  & 51 & 64 \cite{M3} & 61\\ 
7  & 93 & 110 \cite{M3} & 105\\
8  & 183 & 186 \cite{M3} & 183 \\
9  & & 309 \cite{M3} & 297\\
10 & & &  472
\end{tabular}
\\[0.3cm]
Table 1.  Bounds on $B(n)$.
\end{table}

In this table, the values in the column of the best lower bounds known
correspond to the number of points in an hemisphere from the best
known kissing configurations, given by the root systems $D_3$, $D_4$,
$D_5$, $E_6$, $E_7$, $E_8$.

Our method gives a tight upper bound in three cases. In dimension $3$
we get with parameters $d = N = 4$ the bound $B(3) \leq 9.6685$ and
hence we recover the exact values $B(3) = 9$ first proved by G. Fejes
T\'oth (\cite{F}). In dimension $4$ we get with parameters $d = N = 6$
the bound $B(4) \leq 18.5085$ and hence we recover the exact value
$B(4)=18$ first proved by O.R. Musin (\cite{M2}). In dimension $8$ we
find a new tight upper bound. The famous configuration of $240$ points
of $S^7$ given by the root system $E_8$ is well known to be an optimal
spherical code of minimal angular distance $\pi/3$, which is moreover
unique up to isometry. Optimality is due to A.M. Odlyzko and
N.J.A. Sloane (\cite{OS}), and independently to V.I. Levenshtein
(\cite{Le}), uniqueness is due to E. Bannai and N.J.A. Sloane
(\cite{BS}).  From these $240$ points we get a code of the hemisphere
as follows: Take $e$ among these points, then the subset of those
points lying in the hemisphere with center $e$ consists in $183$
points. We obtain a bound of $183.012$ with $d=N=8$ in our
computation. Hence, it proves that it is a maximal code of the
hemisphere, in other words that
\[
B(8)=183.
\]
It is reasonable to believe that the configuration of $183$ points
of $E_8$ is unique up to isometry. Unfortunately we were not able to prove it.

\section{Polynomials}
\label{polynomials}

\subsection{Polynomial restatement of the SDP bound for codes in spherical caps.}

We want to give an equivalent expression of the bound provided by
Theorem \ref{Th SDP} in terms of polynomials. Such an expression will
be useful to prove analytic bounds without the use of software for
solving semidefinite programs, just like in the case of the linear
programming (LP) bound (see e.g.~\cite{OS}). Moreover, we aim at
setting bounds in the form of explicit functions of $\cos\theta$ and
$\cos\phi$. We start with a lemma which shows that any polynomial in
the variables $u,v,t$ can be expressed in terms of the matrix
coefficients of the $Y_k^n(u,v,t)$.  In our situation it suffices to
restrict to polynomials which are symmetric in $u,v$. We introduce the
following notation:
\[
R_d:=\{F\in \R[u,v,t] : F(u,v,t)=F(v,u,t), \deg_{(u,t)}(F)\leq  d, \deg_{t}(F) \leq d\},
\]
where $\deg_{(u,t)}$ stands for the total degree in the variables
$u,t$.

\begin{lemma} Let $F(u,v,t)\in R_d$. There exists a unique 
sequence of $d+1$ real
symmetric matrices $(F_0,F_1,\dots, F_d)$ such that
$F_k$ is a  $(d-k+1)\times  (d-k+1)$ matrix and
\begin{equation}\label{dec}
F(u,v,t)=\sum_{k=0}^d \langle F_k , \Y(u,v,t)\rangle.
\end{equation}
We shall say that $(F_0,\dots,F_d)$ are {\em the matrix coefficients} of $F$.
\end{lemma}

\begin{proof} 
The polynomials $Q^{n-1}_k(u,v,t)$ have degree $k$ in the variable
$t$. Hence, $F(u,v,t)$ has a unique expression of the form
\[
F(u,v,t) = \sum_{k=0}^d q_k(u,v) Q^{n-1}_k(u,v,t),
\]
where $q_k(u,v)$ is symmetric in $u,v$ and has degree in $u$ at most
$d-k$. Since $P_i^{n+2k}(u)$ has degree $i$, $q_k$ has a unique
expression as a linear combination of the products $\lambda_{i,j}
P_i^{n+2k}(u)P_j^{n+2k}(v)$ for $0\leq i,j\leq d-k$. Thus, there is a
symmetric $(d-k+1) \times (d-k+1)$ matrix $F_k$ so that
\[
q_k(u,v) = \sum_{0 \leq i,j \leq d-k} (F_k)_{i,j} \lambda_{i,j}
P_i^{n+2k}(u)P_j^{n+2k}(v).
\]
Since one can write $Y_k(u,v,t)$ as $Q^{n-1}_k(u,v,t) (\lambda_{i,j}
P_i^{n+2k}(u)P_j^{n+2k}(v))$ we obtain decomposition~\eqref{dec}.
\end{proof}

\begin{remark} 
The matrix coefficients of a polynomial $F$ do only trivially depend
on the choice of $d$.  The matrix coefficients associated to $d' \geq
d$ will simply be the ones associated to $d$, enlarged by sufficiently
many rows and columns of zeros.
\end{remark}

\begin{remark}
\label{rq diag} 
From \cite[Proposition 3.5]{BV}, the polynomials $P_k^n(t)$ are linear
combinations of diagonal elements of the matrices $\Y$ with non
negative coefficients. As a consequence, the matrix coefficients of
any polynomial $P(t)\in \R[t]$, are diagonal matrices.  If $P(t)=\sum
f_k P_k^n(t)$, with all $f_k \geq 0$, then the matrix coefficients
$F_k$ of $P$ are also non negative, and, moreover, the top left corner
of $F_0$ equals $f_0$.
\end{remark}

The following reformulation of Theorem \ref{Th SDP} is an analogue of
the classical expression of the linear programming bound (see
e.g.~\cite[Chapter 9, Theorem 4]{CS}).

\begin{theorem}
\label{pol bound}
Let $E_0$ be the matrix whose only non zero entry is the top left
corner which contains $1$.  For a polynomial $F(u,v,t)\in R_d$ let
$(F_0,\dots,F_d)$ be symmetric matrices such that
\[
F(u,v,t)=\sum_{k=0}^d \langle F_k , \Y(u,v,t)\rangle.
\]
Suppose the following conditions hold:
\begin{enumerate}
\item[(a)] $F_k\succeq 0$ for all $0\leq k\leq d$.
\item[(b)] $F_0-f_0E_0 \succeq 0$ for some $f_0>0$.
\item[(c)] $F(u,v,t)\leq 0$ for all $(u,v,t)\in \Delta$.
\item[(d)] $F(u,u,1)\leq B$ for all $u\in [\cos\phi, 1]$.
\end{enumerate}
Then, for any code $C$ in $\CP$ with minimal angular distance at least $\theta$,
\[
\card(C)\leq \frac{B}{f_0}.
\]
\end{theorem}

\begin{proof} 
The statement follows immediately from Theorem \ref{Th SDP} because
the matrices $G_0=F_0/f_0-E_0$ and $G_k=F_k/f_0$ for $1\leq k\leq d$
are a feasible solution to the SDP \eqref{dual sdp} with $M = B/f_0 -
1$.

We also give a direct proof, which has the additional feature to give
information about the case when the obtained bound coincides with the
size of a certain code. Let
\[
S:=\sum_{(c,c')\in C^2} F(e\cdot c, e\cdot c',c\cdot c').
\]
We expand $F$ in the $\Y$'s:
\[
S=\sum_{k=0}^d \langle F_k, \sum_{(c,c')\in C^2} 
\Y(e\cdot c, e\cdot c',c\cdot c')\rangle.
\]
On one hand, from the positivity property \eqref{pos Y} together with the fact that
$\langle A, B \rangle \geq 0$ for any two positive semidefinite matrices
$A,B$ we obtain
\begin{equation}
\label{1st ineq}
\begin{array}{lll}
 S & \geq & \langle f_0E_0, \displaystyle\sum_{(c,c')\in C^2} \overline{Y}_0^n(e\cdot c,
e\cdot c',c\cdot c')\rangle\\
& = & f_0  \displaystyle \sum_{(c,c')\in C^2} \big(\overline{Y}_0^n\big)_{0,0}(e\cdot c,
e\cdot c',c\cdot c') \;=\; f_0\card(C)^2.
\end{array}
\end{equation}
On the other hand, if we split the sum $S$ into diagonal terms
belonging to pairs $(c,c)$ and into cross terms belonging to pairs
$(c,c')$ with $c \neq c'$, we obtain from condition (c) and (d)
\begin{equation}
\label{2nd ineq}
\begin{array}{llclc}
S & = & \displaystyle \sum_{c\in C} F(e\cdot c, e\cdot c,1) & + &
\sum_{(c,c')\in C^2, c\neq c'} F(e\cdot c, e\cdot c',c\cdot c')\\[1ex]
& \leq & B \card(C) & + & 0,
\end{array}
\end{equation}
because $(e\cdot c, e\cdot c,1)\in \Delta_0$ and $(e\cdot c, e\cdot
c',c\cdot c')\in \Delta$ if $c\neq c'$. Now \eqref{1st ineq} and
\eqref{2nd ineq} together give the inequality $\card(C)\leq B/f_0$.
\end{proof}

\begin{remark} 
Like in the LP method, the above proof gives additional information on
the case of equality. Namely, if for a given code $C$ and a given
polynomial $F$, we have $\card(C) = B/f_0$, the inequality \eqref{2nd
ineq} must be an equality. So, $F(u,v,t)=0$ for all $(u,v,t)$ running
through the set of triples $(e\cdot c,e\cdot c', c \cdot c')$ with
$c\neq c'$ and $(c,c')\in C^2$, and $F(u,u,1)=B$ for all $u=e\cdot c$
with $c\in C$.
\end{remark}

\begin{remark} 
In view of explicit computations, it is more convenient to remove the
factor $\lambda_{i,j}$ from the coefficients of $Y_k^n$, so that
polynomials with rational coefficients have rational matrix
coefficients. It changes the above defined $F_k$ to congruence, hence
does not affect the property to be positive semidefinite. These are
the matrix coefficients we discuss about in the next two examples.
\end{remark}

\noindent 
{\bf Example 1.} ($d=1$)

\noindent
We consider the polynomial $F=t-\cos\theta-uv+\cos^2\phi$. The
matrices of the decomposition \eqref{dec} are:
$F_0=\left(\begin{smallmatrix} a& 0\\0&0\end{smallmatrix}\right)$ with
$a=cos^2\phi-\cos\theta$ and $F_1=\left(\begin{smallmatrix}
1\end{smallmatrix}\right)$. Condition (a) of Theorem \ref{pol bound}
is fulfilled if $a\geq 0$. Condition (b) holds for $f_0=a$. Obviously
(c) holds if $\cos\phi\geq 0$ and $B=1-\cos\theta$ because
$F(u,u,1)=1-\cos\theta-u^2+\cos^2\phi$. We obtain:
\[
\text{If }\cos\phi\geq 0\text{ and  }\cos\theta< \cos^2\phi, \text{ then }
A(n,\theta,\phi)\leq \frac{1-\cos\theta}{\cos^2\phi-\cos\theta}.
\]

It is worth to point out that the polynomial 
$G=(t-\cos\theta)-\cos\phi(u+v-2\cos\phi)$ leads to
exactly the same bound. This time 
$F_0=\left(\begin{smallmatrix} c+a& -c\\-c&1\end{smallmatrix}\right)$ 
with $c=\cos\phi$, $f_0=a$, $B=1-\cos\theta$. 

The above bound is already proved in \cite[Theorem 5.2]{BM}. Indeed
with the notations of \cite{BM}, let $w(\theta,\phi)$ be defined by
$\cos w(\theta, \phi)=(\cos\theta-\cos^2\phi)/(\sin^2\phi)$; we have
just proved that the Rankin bound for $A(n-1,w(\theta,\phi))$ is also
a bound for $A(n,\theta,\phi)$. More generally, LP bounds for
$A(n-1,w(\theta,\phi))$ are also bounds for $A(n,\theta,\phi)$: Let $f(x)$ be a
polynomial of degree $d$ that realizes an LP bound on $S^{n-2}$
for the angle $w(\theta,\phi)$. We can take polynomial approximations
of the function
\[
F(u,v,t)=\big((1-u^2)(1-v^2)\big)^{d/2}f\Big(\frac{t-uv}{\big((1-u^2)(1-v^2)\big)^{1/2}}\Big)
\]
obtained by the truncated developments of the powers
$\big((1-u^2)(1-v^2)\big)^{k/2}$ around $u=\cos\phi$, $v=\cos\phi$.

\bigskip

\noindent 
{\bf Example 2.} ($d=2$)
\\
\noindent
We consider the polynomial $F=(t+1)(t-\cos\theta) +
a\big((u-\cos\phi)(u-1)+(v-\cos\phi)(v-1)\big)$. The parameter $a>0$
will be chosen later to optimize the bound. Condition (c) is obviously
fulfilled and condition (d) holds with $B=2(1-\cos\theta)$. The
polynomial $(t+1)(t-\cos\theta)$ has non negative coefficients if we
expand it in terms of the basis $P^n_k(t)$ whenever $\cos\theta \leq
1/n$. More precisely its constant coefficient
equals $\big(\frac{1}{n}-\cos\theta\big)$ while the two others are
positive. So we only need to make sure that $F_0$ is positive
semidefinite.  We find that:
\[
F_0=\left( \begin{matrix}
2a(\frac{1}{n}+\cos\phi)+\frac{1}{n}-\cos\theta & -a(1+\cos\phi) &
a(1-\frac{1}{n})\\
-a(1+\cos\phi) & (1-\cos\theta) & 0\\
a(1-\frac{1}{n}) & 0 &(1-\frac{1}{n})
\end{matrix}
\right).\]
Let
\[
f_0(a):=-a^2\Big(\frac{(1+\cos\phi)^2}{1-\cos\theta}+(1-\frac{1}{n})\Big)+2a\big(\frac{1}{n}+\cos\phi\big)+\frac{1}{n}-\cos\theta.
\]
Then, an easy calculation shows that $F_0\succeq 0$ iff $f_0(a)\geq
0$, and that $F_0-f_0E_0\succeq 0$ iff $f_0\leq f_0(a)$.
The best bound is obtained when $f_0=f_0(a)$ attains the maximal value
\[
(f_0)_{\max}=\big(\frac{1}{n}-\cos\theta\big) +
\frac{\big(\frac{1}{n}+\cos\phi\big)^2}
{\Big(\frac{(1+\cos\phi)^2}{1-\cos\theta}+1-\frac{1}{n}\Big)}.
\]
The final bound equals
\[
\frac{2(1-\cos\theta)}{(f_0)_{\max}}.
\]
and is valid as long as $(f_0)_{\max}>0$ and
$\big(\frac{1}{n}+\cos\phi\big)>0$ (this last condition holds because
$(f_0)_{\max}$ must be attained at a positive $a$).

It is worth noticing that the resulting bound is smaller than the LP
bound for the entire sphere $A(n,\theta)$, obtained from the
polynomial $(t+1)(t-\cos\theta)$, which is
\[
\frac{2(1-\cos\theta)}{\big(\frac{1}{n}-\cos\theta\big)}.
\]
For example, when $\cos\phi=\cos\theta=0$, we recover the exact bound
of $2n-1$ (see also \cite{Ku}).  

\begin{remark}
We can interpret the two examples treated above as follows: in both
cases, we have perturbed the optimal polynomial for the LP method,
respectively $t-\cos\theta$ and $(t+1)(t-\cos\theta)$, by a
polynomial in the variables $u,v$, which affects the first matrix
coefficient $F_0$ and increases the value of the constant coefficient
$f_0$.  However it seems difficult to generalize this approach.
\end{remark}

\subsection{Orthogonality relations.}

In this subsection, we calculate the scalar product induced on
$\R[u,v,t]$  by the natural scalar product on $\Pl(S^{n-1})$ defined
by \eqref{ip}.

\begin{proposition}\label{ch var}
Let $P\in \R[u,v,t]$ be a polynomial. We have 
\begin{align*}
&\frac{1}{\omega_n^2}\int_{(S^{n-1})^2} P(e\cdot x, e\cdot y, x\cdot y)
  d\omega_n(x) d\omega_n(y)=\int_{\Omega} P(u,v,t) k(u,v,t) du dv dt
\end{align*}
where
\begin{equation*}
k(u,v,t)= \frac{\omega_{n-1}\omega_{n-2}}{\omega_n^2}(1-u^2-v^2-t^2+2uvt)^{\frac{n-4}{2}}
\end{equation*}
and
\[
\begin{array} {ll}
\Omega=\{(u,v,t) : &-1 \leq u, v,t\leq 1, \\[1ex]
                   &1+2uvt-u^2-v^2-t^2 \geq 0\}.
\end{array}
\]
\end{proposition}

\begin{proof} If $u=e\cdot x$ and $\zeta\in S^{n-2}$ is defined by
$x=ue+(1-u^2)^{\frac{1}{2}}\zeta$, we have
\[
d\omega_n(x)=(1-u^2)^{\frac{n-3}{2}}dud\omega_{n-1}(\zeta).
\]
With $y=ve+(1-v^2)^{\frac{1}{2}}\xi$, we have
\begin{align*}
&\int_{S^{n-1}} P(e\cdot x, e\cdot y, x\cdot y)d\omega_n(x) \\
&=\int_{S^{n-2}}\int_{-1}^1
P(u,v,uv+\big((1-u^2)(1-v^2)\big)^{\frac{1}{2}}\zeta\cdot \xi)(1-u^2)^{\frac{n-3}{2}}dud\omega_{n-1}(\zeta)\\
&= \omega_{n-2}\int_{-1}^1 \int_{-1}^1
P(u,v,t)(1-\alpha^2)^{\frac{n-4}{2}}(1-u^2)^{\frac{n-3}{2}}d\alpha du,
\end{align*}
where $t:=uv+\big((1-u^2)(1-v^2)\big)^{\frac{1}{2}}\alpha$. With this
change of variables having Jacobian
$\big((1-u^2)(1-v^2)\big)^{\frac{1}{2}}$ we obtain
\begin{eqnarray*}
& &\int_{S^{n-1}} P(e\cdot x, e\cdot y, x\cdot y)d\omega_n(x) \\
&=& \omega_{n-2}\int_{\Omega(v)}
P(u,v,t) (1-u^2-v^2-t^2+2uvt)^{\frac{n-4}{2}}(1-v^2)^{-\frac{n-3}{2}}du dt,
\end{eqnarray*}
where
\[
\begin{array} {ll}
\Omega(v)=\{(u,t) : &-1 \leq u,t\leq 1, \\[1ex]
                   &1+2uvt-u^2-v^2-t^2 \geq 0\}.
\end{array}
\]
Hence
\begin{eqnarray*}
& &\int_{(S^{n-1})^2} P(e\cdot x, e\cdot y, x\cdot y)d\omega_n(x)d\omega_n(y) \\
&=&\omega_{n-1}\omega_{n-2}\int_{\Omega}
P(u,v,t) (1-u^2-v^2-t^2+2uvt)^{\frac{n-4}{2}}du dv dt
\end{eqnarray*}
\end{proof}

\begin{definition} With the notations of Proposition \ref{ch var},
the following expression defines a scalar product
on $\R[u,v,t]$:
\begin{equation}\label{def ip}
[F,G]=\int_{\Omega} F(u,v,t)G(u,v,t)k(u,v,t) du dv dt.
\end{equation}
From Proposition \ref{ch var}, it is the scalar product induced 
by the standard scalar product
\eqref{ip}
on $\Pl(S^{n-1})$.
\end{definition}

The subspaces $H_{k,i}^{n-1}$ are pairwise orthogonal. Consequently
the matrix coefficients of $Y_k^n(u,v,t)$ are pairwise orthogonal for
$[\cdot,\cdot]$. Their norm is also easy to compute, and we obtain the
following useful formulas:

\begin{proposition}
\noindent
\begin{enumerate}
\item[(a)] For all $k, k'$ and all $i, j$, $i', j'$ we have 
\begin{equation}\label{ip Y}
[\big(Y_k^n\big)_{i,j}, \big(Y_{k'}^n\big)_{i',j'}]=
\frac{\delta_{(i,j,k), (i',j',k')}}{h_k^{n-1}}.
\end{equation}

\item[(b)] For all symmetric matrices $A,B$ and all $k, k'$ we have
\begin{equation}\label{ip tr}
[\langle A, \Y \rangle, \langle B, \overline{Y}_{k'}^n \rangle]= 
\frac{\delta_{k,k'}\langle A, B\rangle  }{h_k^{n-1}}.
\end{equation}
\end{enumerate}
\end{proposition}

\begin{proof} 
Obvious. 
\end{proof}

\subsection{Characterization of the positive definite polynomials.}

In view of Theorem \ref{pol bound}, we are concerned with the
construction of polynomials satisfying condition (a). We prove in this
subsection that this property is stable under multiplication. We start
with a characterization of the set of polynomials satisfying (a)  of Theorem
\ref{pol bound}.

\begin{definition} We say that the polynomial $F(u,v,t)\in \R[u,v,t]$
is {\em positive definite} if, for all finite $C\subset S^{n-1}$,
for all functions $\alpha:C\to \R$,
\begin{equation}\label{pos prop}
\sum_{(c,c')\in C^2} \alpha(c)\alpha(c') F(e\cdot c, e\cdot c', c\cdot
c')\geq 0.
\end{equation}
\end{definition}

The polynomials $F(u,v,t)$ of the form
\[
F(u,v,t)=\sum_{k=0}^d \langle F_k , \Y(u,v,t)\rangle
\]
with $F_k\succeq 0$ for all $0\leq k\leq d$ are positive definite in the
above sense. Note that \eqref{pos prop} is slightly stronger than the positivity
property of the matrices $Y_k^n$ proved in Theorem \ref{Th pos Y}; the 
argument is essentially the same, as it follows from the equality
\[
\displaystyle \sum_{(c,c')\in C^2}\alpha(c)\alpha(c')Z_k^n(c,c')=\Big(\sum_{c\in
C}\alpha(c)E_k^n(c)\Big)\Big(\sum_{c\in
C}\alpha(c)E_k^n(c)\Big)^t\succeq 0.
\] 
We prove with next proposition that all positive definite polynomials in
$R_d$ arise in this way. 

\begin{proposition}\label{Prop pos prop}
Let $F(u,v,t)\in R_d$.  Let $(F_0,\dots,F_d)$ be symmetric matrices
such that
\[
F(u,v,t)=\sum_{k=0}^d \langle F_k , \Y(u,v,t)\rangle.
\]
If $F$ is positive definite, then $F_k\succeq 0$ for all $0\leq k\leq
d$.

\end{proposition}

\begin{proof} Let $\tilde F(x,y)=F(e\cdot x, e\cdot y,
  x\cdot y)$. By compactness, $F$ is positive definite if and only if
  for all $f\in \Pl(S^{n-1})$,
\[
\int_{\big(S^{n-1}\big)^2} f(x)f(y) \tilde F(x,y) d\omega_n(x) d\omega_n(y) \geq 0.\]
As a consequence, if $Q(x)$ is any matrix,
\[
\int_{\big(S^{n-1}\big)^2} \langle Q(x), Q(y)\rangle \tilde F(x,y)
d\omega_n(x) d\omega_n(y)\geq 0.\] Let us fix $k \in \{0, \ldots, d\}$
and let $A$ be a $(d-k+1)\times (d-k+1)$ symmetric, positive
semidefinite matrix. Because of expression \eqref{Z} of $Z_k^n$, we
can write $\langle A, Z_k^n(x,y)^t\rangle$ in the form $\langle Q(x), Q(y)\rangle$. Hence,
\[
\int_{\big(S^{n-1}\big)^2} \langle A, Z_k^n(x,y)^t\rangle \tilde F(x,y) d\omega_n(x) d\omega_n(y)\geq 0.\]
In terms of the scalar product $[\cdot,\cdot]$ this is equivalent to
\[
[\langle A, \Y\rangle, F]\geq 0.
\]
Since from \eqref{ip
  tr} $[\langle A, \Y \rangle, F]=\big(h_k^{n-1}\big)^{-1}\langle A, F_k\rangle$,
we have proved that $ \langle A, F_k\rangle\geq 0$ for all $A \succeq
  0$, and so $F_k\succeq 0$.
\end{proof}

\begin{remark} This characterization of positive definite functions is
  in fact already proved in \cite[Section III]{Boch} in a more general
  context: for compact spaces which are homogeneous under the action
  of their automorphism group, but not necessarily two-point
  homogeneous. The assumption that the group acts transitively is
  however not needed in the proof.
\end{remark}

\begin{corollary}\label{prod pos} Let $F, G\in R_d$.
If $F$ and $G$ are positive definite, then the 
product $FG$ is also positive definite.
\end{corollary}

\begin{proof} From Proposition \ref{Prop pos prop} it suffices to consider
  the case $F=\langle A, \Y \rangle$, $G=\langle B, \overline{Y}_l^n
  \rangle$, where $A$ and $B$ are positive semidefinite matrices.
  Again, we write $\langle A, Z_k^n(x,y)^t\rangle=\langle Q(x), Q(y)\rangle$ and
  $\langle B, Z_l^n(x,y)^t\rangle=\langle T(x), T(y)\rangle$.  With the formula
\[
\langle Q(x), Q(y) \rangle \langle T(x),T(y)\rangle=\langle Q(x)\otimes
T(x), Q(y)\otimes T(y)\rangle
\]
we have

\begin{align*}
\sum_{(c,c')\in C^2} \alpha(c)\alpha(c')&\tilde F(c,c')\tilde G(c,c')\\
&=\sum_{(c,c')\in C^2 }\alpha(c)\alpha(c')\langle Q(c), Q(c')\rangle
  \langle T(c), T(c')\rangle \\
&=\sum_{(c,c')\in C^2} \langle \alpha(c)Q(c)\otimes
T(c), \alpha(c')Q(c')\otimes T(c'))\rangle \\
&=\langle U_C, U_C \rangle\geq 0
\end{align*}
with
\[
U_C=\sum_{c\in C} \alpha(c)Q(c)\otimes T(c).
\]

\end{proof}

\subsection{Reproducing kernels.}

We define the kernel $K_d : \R^3 \times \R^3 \to \R$ by
\begin{equation}\label{def k}
K_d((u,v,t),(u',v',t')):= \sum_{k=0}^n h_k^{n-1} \langle \Y(u,v,t), \Y(u',v',t') \rangle.
\end{equation}

\begin{proposition}
The kernel $K_d$ is the {\em reproducing kernel} of the space 
$R_d$, i.e., for all $F\in R_d$ and all $(u',v',t') \in \R^3$ we have
\begin{equation}\label{rep k}
[K_d(\cdot,(u',v',t')), F] = F(u',v',t').
\end{equation}
\end{proposition}

\begin{proof} It is straightforward from \eqref{ip tr}.

\end{proof}


\begin{thebibliography}{[15]}

\bibitem{B} C. Bachoc,
{\em Linear programming bounds for codes in Grassmannian spaces},
IEEE Trans. Inf. Th.  {\bf 52-5} (2006), 2111--2125.

\bibitem{BV} C. Bachoc, F. Vallentin,
{\em New upper bounds for kissing numbers from semidefinite
  programming}, to appear in J. Amer. Math. Soc.

\bibitem{BM} A. Barg, O.R. Musin, {\em  Codes in spherical caps},
Adv. Math. Comm. {\bf 1} (2007), 131--149.

\bibitem{Boch}
S. Bochner,
{\em Hilbert distances and positive definite functions},
Ann. of Math. {\bf 42} (1941) 647--656.

\bibitem{Bo}
B. Borchers,
{\em CSDP, A C library for semidefinite programming},
Optimization Methods and Software {\bf 11} (1999) 613--623.

\bibitem{BS} E. Bannai, N.J.A. Sloane,
{\em Uniqueness of certain spherical codes},
Canad J. Math. {\bf 33} (1981), 437--449.

\bibitem{CS} J.H. Conway, N.J.A. Sloane, 
{\em Sphere Packings, Lattices and Groups}, 
Springer-Verlag, 1988.

\bibitem{D} P. Delsarte, 
{\em An algebraic approach to the association schemes of coding
  theory},
Philips Res. Rep. Suppl. (1973), vi+97.

\bibitem{DGS} P. Delsarte, J.M. Goethals, J.J. Seidel,
{\em Spherical codes and designs}, Geom. Dedicata {\bf 6} (1977), 363--388.

\bibitem{F} G. Fejes T\'oth, {\em Ten-neighbor packing of equal balls},
Periodica Math. Hungar. {\bf 12}, (1981) 125--127.


\bibitem{KL} G.A.~Kabatiansky, V.I.~Levenshtein,
{\em Bounds for packings on a sphere and in space},
Problems of Information Transmission {\bf 14} (1978), 1--17.


\bibitem{Ku} W. Kuperberg, {\em Optimal arrangements in packing
congruent balls in a spherical container}, Discrete Comput. Geom. {\bf 37}
(2007), 205--212.

\bibitem{Le} V.I.~Levenshtein,
{\em On bounds for packing in $n$-dimensional Euclidean space},
Soviet Math. Dokl. {\bf 20} (1979), 417--421.

\bibitem{M2} O.R. Musin, {\em The one-sided kissing number in four dimensions},
Periodica Math. Hungar. {\bf 53}, (2006) 209--225.

\bibitem{M3} O.R. Musin, {\em Bounds for codes by
semidefinite programming}, preprint, September 2006,
\href{http://arxiv.org/abs/math.CO/0609155}{arXiv:math.MG/0609155}.

\bibitem{OS} A.M. Odlyzko, N.J.A. Sloane,
{\em New bounds on the number of unit spheres that can touch a unit
  sphere in n dimensions}, 
J. Combin. Theory Ser. A {\bf 26} (1979), 210--214.

\bibitem{Pu} M. Putinar,
{\em Positive polynomials on compact semi-algebraic sets},
Ind. Univ. Math. J. {\bf 42} (1993), 969--984.
\end{thebibliography}
\end{document}